\magnification=\magstep1
\documentstyle{amsppt}
\UseAMSsymbols
\voffset=-3pc
\loadbold
\loadmsbm
\loadeufm
\parskip=6pt
\def\pop{\partial\overline\partial}
\def\bC{\Bbb C}
\def\bR{\Bbb R}
\def\cH{\Cal H}
\def\cO{\Cal O}
\def\Im{\text{Im}\,}
\def\Re{\text{Re}\,}
\def\id{\text{id}}
\def\Ker{\text{Ker }}
\NoBlackBoxes
\topmatter
\title Weak Geodesics in the Space of K\"ahler Metrics
\endtitle
\author
Tam\'as Darvas and L\'aszl\'o Lempert\footnote""{$^*$Research supported by NSF grant DMS0700281.
Part of the research was done while the second author enjoyed the hospitality of the
Universit\'e Pierre et Marie Curie, Paris.\hfill\break}\endauthor
\leftheadtext{Space of K\"ahler metrics}
\rightheadtext{Tam\'as Darvas and L\'aszl\'o Lempert}
\address
Department of Mathematics, Purdue University, West Lafayette, IN\ 47907
\endaddress
\abstract
Given a compact K\"ahler manifold $(X,\omega_0)$, according to Mabuchi, the set $\cH_0$ of K\"ahler forms
cohomologous to $\omega_0$ has the natural structure of an infinite dimensional Riemannian manifold.
We address the question whether points in $\cH_0$ can be joined by a
geodesic, and strengthening the
finding of \cite{LV}, we show that this cannot always be done even with a
certain type of generalized geodesics.
As in \cite{LV}, the result is obtained through the analysis of a
Monge--Amp\`ere equation.
\endabstract
\endtopmatter
\TagsOnLeft
\document
\subhead 1.\ Introduction\endsubhead

Let $X$ be a connected compact complex manifold of dimension $m>0$ and
$\omega_0$ a smooth K\"ahler form on it.
In the 1980s Mabuchi discovered that there is a natural infinite dimensional Riemannian manifold structure on the
set $\cH_0$ of smooth K\"ahler forms cohomologous to $\omega_0$, and on the set
$$
\cH=\{v\in C^\infty(X) : \omega_0+i\pop v>0\}
$$
of smooth strongly $\omega_0$--plurisubharmonic functions.
He also showed that $\cH$ is isometric to the Riemannian product $\cH_0\times\bR$, \cite{M}.
In \cite{LV}, answering a question posed by Donaldson, Vivas and the second author proved that in general there is no geodesic of class $C^2$ between two points in $\cH$, resp. in $\cH_0$; in fact, there is not even
one of Sobolev regularity $W^{1,2}$.

Since geodesics and their generalizations, weak geodesics, potentially play an important role in the study of
special K\"ahler metrics (for geodesics, see \cite{D1, M}), it is of interest to know whether two points in $\cH$ can
be connected at least by a weak geodesic.
What the notion of weak geodesic should be is suggested by Semmes' reformulation of the geodesic equation in $\cH$, see
\cite{S}.
Let $S=\{s\in\bC: 0 < \Im s < 1\}$ and $\omega$ the pullback of $\omega_0$ by the projection $\overline S\times
X\to X$.
With any $C^2$ curve $[0,1]\ni t\mapsto v_t\in\cH$ associate a function $u:\overline S\times X\to\bR$,
$$
u(s,x)=v_{\,\Im s}(x),
$$
itself a $C^2$ function.
Then $t\mapsto v_t$ is a geodesic if and only if $u$ satisfies the Monge--Amp\`ere equation $(\omega+i\pop u)^{m+1}=0$.
Therefore a $C^2$ geodesic connecting $0, v\in\cH$ gives rise to a solution
$u\in C^2(\overline S\times X)$ of a
boundary value problem for this Monge--Amp\`ere equation on
$\overline S\times X$; furthermore $\omega+i\pop u\geq 0$.
This latter is expressed by saying that $u$ is $\omega$--plurisubharmonic.
By a weak, or generalized, geodesic connecting, say, $0,v\in\cH$ one then means an $\omega$--plurisubharmonic solution
$u:\overline S\times X\to\bR$ of the problem
$$
\aligned
(\omega+i\pop u)^{m+1}&=0,\\
u(s+\sigma,x)&=u(s,x),\quad\text{ if }(s,x)\in\overline S\times X,\ \sigma\in\bR,\\
u(s,x)&=\cases 0,&\text{if }\,\Im s=0\\ v(x),&\text{if }\,\Im s=1.\endcases
\endaligned\tag1.1
$$
It has to be assumed that $u$ is sufficiently regular so that $(\omega+i\pop u)^{m+1}$ can be given sense; for example,
according to \cite{BT}, the continuity of $u$ more than suffices.
X.X.~Chen has indeed proved that for $v\in\cH$ (1.1) admits a continuous $\omega$--plurisubharmonic solution for which the
current $\pop u$ is represented by a bounded form, see \cite{C} and complements in \cite{B\l}.
In other words, any two points in $\cH$ can be connected by a weak geodesic.
One should keep in mind, though, that a weak geodesic $u$ need not give rise
to a curve in $\cH$, first because $v_t=u(t,\cdot)$
is not necessarily $C^\infty$, not even $C^2$, and second because even if
$v_t$ is $C^\infty$, there is no reason why it should be strongly
$\omega_0$--plurisubharmonic.

In this paper we show that the regularity that Chen obtains cannot be improved:\ (1.1) may have a solution with $\pop u$
bounded, but in general it will not have a solution with $\pop u$ continuous.

If $\overline Z$ is a complex manifold, possibly with boundary, and
$Z=\text {int}\ \overline Z$, we define
$$
\multline
C^{\pop}(\overline Z)=\\
\{w\in C(\overline Z)\colon\text{ the current } \pop(w|Z)
\text{ is represented by a form continuous on } \overline Z\}.
\endmultline
$$
Given $w\in C^{\pop}(\overline Z)$, we will simply write $\pop w$ for the continuous form on $\overline Z$ that represents the current $\pop (w|Z)$, and if
$z_1,z_2,\ldots$ are local coordinates on $Z$, we write $w_{z_j\bar z_k}$ for the coefficient of $dz_j\wedge d\bar z_k$ in $\pop w$.

Clearly $C^2(\overline Z)\subset C^{\pop}(\overline Z)$, and it is well
understood in harmonic analysis that the inclusion is strict. For
example, if $\overline Z=\{\zeta\in\bC:|\zeta|\le 1/2\}$ and $k=2,3,\ldots$, the
function
$$
w(\zeta)=\cases\zeta^k\log\log|\zeta|^{-2}, &\text{if }0<|\zeta|\le 1/2 \\
0,&\text{if }\zeta=0
\endcases
$$
is not in $C^k(\overline Z)$, but $w_{\bar\zeta}\in C^{k-1}(\overline Z)$ and
$w_{\zeta\bar\zeta}\in C^{k-2}(\overline Z)$.

\proclaim{Theorem 1.1}Suppose a connected compact K\"ahler manifold $(X,\omega_0)$ admits a holomorphic isometry $g: X\to X$
with an isolated fixed point, and $g^2=\id_X$.
Then there is a $v\in\cH$ for which (1.1) has no $\omega$--plurisubharmonic solution $u\in C^{\pop}(\overline S\times X)$.
One can choose $v$ to satisfy $g^* v=v$.
\endproclaim

The proof will show that among symmetric potentials the $v\in\cH$ in Theorem 1.1 even form an open set.

Theorem 1.1 corresponds to \cite{LV, Theorem 1.2}, but the $C^3$ regularity from \cite{LV} has been lowered.
The proofs here and in \cite{LV} are similar in that, denoting by $x_0\in X$ an isolated fixed point of $g$, in both proofs we
analyze the behavior of a regular solution $u$ in a neighborhood of $\overline S\times \{x_0\}$.
The upshot of the analysis is a condition on the Hessian of the boundary value at $x_0$, a condition that not all $v\in\cH$
satisfy.
In \cite{LV} the analysis involved the Monge--Amp\`ere foliation associated
with a $u\in C^3(\overline S\times X)$, and it
was crucial that the foliation was of class $C^1$.
The foliation method is not available when $u$ is only $C^{\pop}$, and we will have to be thriftier with our tools, but in spite of this, we will recover the same condition on the Hessian as in \cite{LV} when $m=1$. When $m > 1$, the present condition is slightly stronger than the one in \cite{LV}.

\subhead 2.\ Generalities\endsubhead

In this section we collect a few simple facts concerning currents and the homogeneous Monge--Amp\`ere equation.

\proclaim{Proposition 2.1}Let $f\colon Y\to Z$ be a holomorphic map of complex manifolds, $\varphi$ and $\psi$ continuous
forms on $Z$ satisfying $\pop\varphi=\psi$ as currents.
Then $\pop f^*\varphi=f^*\psi$ as currents.
\endproclaim

\demo{Proof}We can assume $Z$ is an open subset of some $\bC^n$.
Regularizing $\varphi$ and $\psi$ by convolutions gives rise to sequences of smooth forms $\varphi_k$ and
$\psi_k=\pop\varphi_k$ that converge locally uniformly to $\varphi$, resp.~$\psi$.
Therefore $f^*\varphi_k\to f^*\varphi$ and $\pop f^*\varphi_k=f^*\pop\varphi_k\to f^*\psi$ locally uniformly, whence the claim
follows from the continuity of $\pop$ in the space of currents.
\enddemo

Next consider a complex manifold $Z$ and a plurisubharmonic $U\in C^{\pop}(Z)$.
Suppose $Y\subset Z$ is a one dimensional, not necessarily closed complex submanifold and $TY\subset \Ker \pop U$.
The normal bundle $NY=(T^{1,0} Z|Y)/T^{1,0} Y$ is a holomorphic vector bundle and $\pop U$ induces a possibly degenerate
Hermitian metric $h$ on it.
With $p\colon T^{1,0} Z|Y\to NY$ the canonical projection, the metric is
$$
h(p\zeta)=\pop U(\zeta,\overline\zeta)\geq 0,\quad \zeta\in T^{1,0} Z|Y.
$$
Thus $h$ is continuous, but can degenerate, i.e., vanish on nonzero vectors
as well.

\proclaim{Proposition 2.2}The metric $h$ is seminegatively curved in the sense that log $h\circ\sigma$ is subharmonic for
any holomorphic section $\sigma$ of $NY$ over some open $Y'\subset Y$.
(Here it is convenient not to exclude from subharmonic functions those that are
identically $-\infty$ on some component of $Y'$.)
Further, on the line bundle $\det NY$ the metric induced by $h$ is
also seminegatively curved.\endproclaim

When $h$ is smooth and nondegenerate, and moreover $(\pop U)^{\dim Z}=0$,
the semi\-negativity of $\det NY$ was first proved by Bedford and Burns in
\cite{BB, Proposition 4.1}, and \cite{CT, Theorem 4.2.8} gives the
seminegativity of $N$ itself.
For possibly degenerate $h$ \cite{BF, Lemma} represents an equivalent result,
albeit without the curvature interpretation, and under the
assumption that $U$ is $C^2$.
Our proof is a variant of the proof in \cite{BF}.

\demo{Proof}For the first statement we only need to prove that
$\log h\circ\sigma$ has the sub--meanvalue property, and this at points where
$h\circ\sigma\neq 0$.
To do so, we can assume $Z\subset\bC^n$ is the unit polydisc,
$Y=Y'=\{z\in Z: z_2=\ldots=z_n=0\}$, and that
$\sigma(z_1,0,0,\ldots)=p(\partial/\partial z_2)$.
Thus
$$
(h\circ\sigma) (z_1,0,\ldots)=U_{z_2\bar z_2} (z_1,0,\ldots).\tag2.1
$$
Green's formula implies for $0<r<1$
$$\multline
\frac 1 {r^2}\ \int_0^1 \big(U(z_1,re^{2\pi it},0,\ldots)-U(z_1,0,0,\ldots)\big)
\, dt
=\\ \frac i{\pi r^2}\ \int_{|z_2|\leq r} (\log r-\log
|z_2|)U_{z_2\bar z_2}(z_1,z_2,0,\ldots)\, dz_2\wedge d\bar z_2,
\endmultline\tag2.2
$$
certainly if $U$ is $C^2$, but then upon regularizing by convolutions, whenever $U$ and $\pop U$ are continuous---as in our
case.
Proposition 2.1, with $f$ the embedding $Y\to Z$, implies $\pop (U|Y)=(\pop U)|Y=0$.
Hence the left hand side of (2.2) is a subharmonic function of $z_1$, and so is the right hand side.
As $r\to 0$, these functions converge locally uniformly to $U_{z_2\bar z_2} (z_1,0,\ldots)$; in light of (2.1)
$h\circ\sigma$ is therefore subharmonic.

If $\varphi\in\cO(Y)$ and $\sigma$ is replaced by $e^{\varphi/2}\sigma$, we obtain that $e^{\Re\varphi}h\circ\sigma$ is also
subharmonic.
Therefore it satisfies the maximum principle, and so does $\Re\varphi+\log h\circ\sigma$; knowing this for all
$\varphi\in\cO(Y)$ is equivalent to the subharmonicity of $\log h\circ\sigma$,
see e.g.~\cite{H, Theorem 1.6.3}.

Now given any holomorphic vector bundle $E\to Y$ of rank $r$, endowed with
a seminegatively curved, possibly degenerate continuous Hermitian metric
$h$, the induced metric on the line bundle $\det E$ is also seminegatively
curved. Indeed, denoting by $h(e,e')$  the inner product of
$e,e'\in E_y$, $y\in Y$,
so that $h(e)=h(e,e)$, for (local) sections $\sigma_1,\ldots\sigma_r$ of $E$
the induced metric is given by
$$
h^{\det}(\sigma_1\wedge\ldots\wedge\sigma_r)=\det\big(h(\sigma_j,\sigma_k)\big).
\tag2.3
$$
If $h$ is smooth and nondegenerate and $y\in Y$, any nonzero holomorphic
section of $\det E$ in a neighborhood of $y$ can be written as $\sigma_1\wedge\ldots\wedge\sigma_r$, where the $\sigma_j$ are holomorphic
sections of $E$ near $y$, and $h(\sigma_j,\sigma_k)$
vanish to second order at $y$ for $j\neq k$. Thus
$\det\big(h(\sigma_j,\sigma_k)\big)-\prod_{j=1}^r h(\sigma_j,\sigma_j)$
vanishes to {\sl fourth} order at $y$.
By virtue of (2.3) this implies that at $y$
$$
i\pop\log h^{\det}(\sigma_1\wedge\ldots\wedge\sigma_r)=
i\pop\log\prod_{j=1}^r h(\sigma_j,\sigma_j)\ge 0.
$$
Therefore $h^{\det}$ is seminegatively curved when $h$ is
smooth and nondegenerate. To prove for a general $h$ we can assume
$Y\subset\bC$ is connected, $E=Y\times\bC^r$ is holomorphically trivial,
and $h^{\det}$ degenerates nowhere. We can regularize $h$ by convolutions,
and obtain
$h^{\det}$ as the locally uniform limit of seminegatively curved metrics,
hence itself seminegatively curved.
\enddemo

Lastly we record a uniqueness result and its corollary:

\proclaim{Proposition 2.3}Given a compact K\"ahler manifold $(X,\omega_0)$ and $v\in\cH$, the equation (1.1) has at most one
$\omega$--plurisubharmonic solution $u\in C^{\pop}(\overline S\times X)$.
\endproclaim

The result follows from \cite{B\l, Proposition 2.2 or Theorem 2.3} or from \cite{PS, p.~144}, once one checks that for
$\omega$--plurisubharmonic $u\in C^{\pop}(\overline S\times X)$ the Monge--Amp\`ere measure $(\omega+i\pop u)^{m+1}$, as
defined e.g.~in \cite{BT}, agrees with what is obtained by taking the exterior power of the continuous form $\omega+i\pop
u$.
Alternatively, the more elementary arguments for \cite{D1, Lemma 6} and the first paragraph of the proof of \cite{LV,
Proposition 2.3} also give uniqueness, provided one first checks the following:\ if $Z$ is a complex manifold and $w\in C^{\pop}(Z)$ is real valued, then
$i\pop w\geq 0$ at any local minimum point of $w$.
Because of Proposition 2.1, it suffices to verify this latter when $\dim Z=1$, and then it is straightforward:\ if $i\pop
w<0$ at a point, then $i\pop w<0$ in a neighborhood, whence $w$ is strongly superharmonic there, and has no local minimum.

\proclaim{Corollary 2.4}Suppose $v\in\cH$ satisfies $g^* v=v$, and $u\in C^{\pop}(\overline S\times X)$ is an
$\omega$--plurisubharmonic solution of (1.1).
Then $u(s,x)=u(s,g(x))$.
\endproclaim

\subhead 3.\ The proof of Theorem 1.1\endsubhead

Let $X,\omega_0,\omega$, and $g$ be as in Theorem 1.1, and let $x_0\in X$ be an isolated fixed point of $g$.
Using \cite{LV, Proposition 2.2} we choose local coordinates $z_1,\ldots,z_m$ in a neighborhood $V\subset X$ of $x_0$ in
which $g$ is expressed as $(z_j)\mapsto (-z_j)$.

\proclaim{Proposition 3.1}If an $\omega$--plurisubharmonic
$u\in C^{\pop}(\overline S\times V)$ solves (1.1) and $u(s,x)=u(s,g(x))$, then
$u(s,x_0)=a\,\Im s$ for $s\in\overline S$, with some $a\in\bR$.
\endproclaim

\demo{Proof \rom{(essentially taken over from \cite{BF, Proposition})}}
The symmetry assumption implies that $u_{s\bar z_j}(s,x_0)=0$, and so
$$
0=(-i\omega+\pop u)^{m+1}=(-i\omega+\partial_X\overline\partial_X u)^m\wedge\ \partial_S\overline\partial_S u\tag3.1
$$
at points of $\overline S\times \{x_0\}$.
Hence for any $s\in\overline S$ either
$(-i\omega+\partial_X\overline\partial_X u)^m$ or
$\partial_S\overline\partial_S u$
vanishes at $(s,x_0)$.
The goal is to show that it is always the latter that vanishes.

We claim that on $S\times \{x_0\}$
$$
\lambda=\log \big(-i\omega+\partial_X\overline\partial_X u\big)^m
\big({\partial\over\partial z_1}\wedge {\partial\over\partial\bar
z_1}\wedge\ldots\wedge\ {\partial\over\partial z_m}\wedge {\partial\over\partial\bar z_m} \big)
$$
is subharmonic and not identically $-\infty$.
Indeed, by Proposition 2.1 $(\pop u)|\{s\}\times X=\pop (u|\{s\}\times X)$ for $s\in S$, and by the continuity of $\pop$,
also for $s\in\overline S$.
But $u(0,\cdot)=0$ is strongly $\omega_0$--plurisubharmonic, hence $\lambda(s,x_0)>-\infty$ when $s=0$, and also when
$s\in\overline S$ is near 0.
As to subharmonicity, it suffices to verify it on the open set
$$
S_0=\{s\in S: \lambda(s,x_0) > -\infty\}.
$$
Choose a smooth $w_0$ in a neighborhood of $x_0\in X$ such that $\omega_0=i\pop w_0$, let $w(s,x)=w_0(x)$ and $U=u+w$.
By what has been observed above,
$U_{s\bar z_j}(s,x_0)=U_{s\bar s} (s,x_0)=0$ if $s\in S_0$; in other
words, $S_0\times \{x_0\}$ is tangential to Ker $\pop U$.
By virtue of Proposition 2.2 $\lambda$ is subharmonic on $S_0\times \{x_0\}$, hence on $S\times \{x_0\}$, as claimed.

Once we know $\lambda|S\times \{x_0\}$ is subharmonic, it follows that $S_0$ is dense in $S$; since by (3.1) $u_{s\bar s}$ vanishes on $S_0\times \{x_0\}$, it vanishes on all of $S\times \{x_0\}$.
The Proposition now follows, because a harmonic function on $S$ that depends only on $\Im s$ must be a linear function of $\Im s$.
\enddemo

In the proof of the next lemma we will make use of the Poisson integral representation of harmonic functions in a strip. If $\psi$ is harmonic $S$, continuous and bounded in $\overline{S}$, then we have the following integral representation (for more on this see \cite{W}):
$$
\psi(\xi + i \eta) = \int_{-\infty}^{+\infty}P(t-\xi, \eta)\psi(t)dt +\int_{-\infty}^{+\infty}P(t-\xi,1- \eta)\psi(t+i)dt,\tag3.2
$$
where $P$ is the following Poisson kernel:
$$P(\xi,\eta)=\frac{\sin  \pi \eta}{2(\cosh\pi \xi - \cos \pi \eta)}.$$
As expected, the above integral representation formula also gives a recipe  to generate bounded continuous harmonic functions in $\overline{S}$ given bounded continuous boundary data on $\partial S$.

Let $\omega_0=\sum^m_{j,k=1}\omega_{jk} dz_j\wedge d\bar z_k$ on $V\subset X$.

\proclaim{Lemma 3.2}Suppose $a\in\bR$ and $u$ is a bounded continuous $\omega$--plurisubharmonic function in $\overline S\times V$ satisfying
$$
\align
u(s,x)&=0,\phantom{(s,x)}\quad \text{ if }(s,x)\in\bR\times V,\\
u(s,x_0)&=a\,\text{Im } s,\quad\, \text{ if }s\in\overline S.
\endalign
$$
If $v=u(i,\cdot)$ is twice differentiable at $x_0$ and $dv=0$ there, then
$$
\big|\sum^m_{j,k=1} v_{z_j z_k} (x_0)\xi_j\xi_k\big|\leq
\sum^m_{j,k=1} \big(2\,\omega_{jk} (x_0)+
v_{z_j\bar z_k}(x_0)\big)\xi_j\overline\xi_k
\qquad\text{for } \xi_j\in\bC.\tag3.3
$$
and this estimate is sharp.
\endproclaim

\demo{Proof}We will assume $a=0$ (otherwise we replace $u(s,x)$ by $u(s,x)-a\,\Im s)$.
Thus $u(s,x_0)=v(x_0)=0$.
By passing to a slice, the proof is reduced to the case $m=1$.
We will denote the local coordinate on $V$ by $z=z_1$; it identifies $V$ and $x_0$ with a neighborhood of $0\in\bC$ and with $0\in\bC$.
Since $m=1$, we need to verify
$$
|v_{zz}(0)|\leq 2\,\omega_{11} (0)+v_{z\bar z}(0).\tag3.4
$$

Suppose $f:\overline{S} \to \Bbb C$ is bounded and holomorphic with $f(  \alpha)=0$ for some $\alpha \in S$. Let $q = v_{zz}(0)$, and choose real numbers $p > v_{z\overline{z}}(0)$ and $r > \omega_{11}(0)$. With a neighborhood  $V'\subset V$ of $0$ we will have
$$
v(z) \leq p |z|^2 + \text{Re } qz^2 \; \text{ and } \;
\omega < i \partial \overline{\partial} r |z|^2
$$
for all $z  \in V'$. 
Clearly, this implies that the function $U(s,z) = r |z|^2 + u(s,z)$ is plurisubharmonic in $S \times V'$ and if $\zeta$ is sufficiently small then
$$\phi(s)=U(s,\zeta f(s)) = r|\zeta f(s)|^2 + u(s,\zeta f(s))$$
is a subharmonic function of  $s \in S$. On the boundary of $\overline{S}$ we have the following estimates:
$$
\phi(s)\ \cases=r|\zeta f(s)|^2,&\text{if $\Im s=0$}\\
\leq (p+r)|\zeta f(s)|^2+\Re q\zeta^2 f(s)^2,&\text{if $\Im s=1$.}\endcases
$$
We take $\zeta$ such that $q\zeta^2$ is nonnegative. Then $\text{Re } q\zeta^2 f(s)^2=|q\zeta^2|\Re f(s)^2$.

Let $\psi_1$,$\psi_2$ and $\psi_3$ be bounded, continuous and harmonic functions on $\overline{S}$ defined by the following boundary data:
$$
\align
\psi_1(s)&=|f(s)|^2, \; \text{ if } s \in \partial S,\\
\psi_2(s)&=0, \;\text{ if Im }s =0 \;\text{ and } \; \psi_2(s)=|f(s)|^2,
\;\text{ if Im }s =1,\\
\psi_3(s)&=0, \;\text{ if Im }s =0\; \text{ and } \;\psi_3(s)=\text{Re }f(s)^2,\; \text{ if Im }s =1.
\endalign
$$
Since $\psi_1$, $\psi_2$, $\psi_3$ and $\phi$ are all bounded on $\overline{S}$, by the maximum principle we obtain 
$$
0 = u(\alpha,0)= \phi( \alpha) \leq r |\zeta|^2\psi_1(\alpha) + p |\zeta|^2 \psi_2(\alpha) + |q||\zeta|^2 \psi_3( \alpha).\tag3.5
$$

We will show that $\psi_2(\alpha)/\psi_1(\alpha)$ can be chosen arbitrarily close to $1/2$ and $\psi_3(\alpha)/\psi_1(\alpha)$ can be chosen arbitrarily close to $-1/2$. We can work with any $\alpha\in S$, but if $\alpha=\xi+i\eta=i/2$, Poisson's formula (3.2) simplifies and gives $\psi_1(i/2)=(I + J)/2$, $\psi_2(i/2)= J/2$ and $\psi_3(i/2)=K/2$, where
$$
\gather I=I(f)=\int_{-\infty}^{+\infty} \frac{|f(t)|^2}{\cosh\pi t }dt,\quad
J=J(f)=\int_{-\infty}^{+\infty} \frac{|f(t + i)|^2}{\cosh\pi t }dt,\\
K=K(f)=\int_{-\infty}^{+\infty} \frac{\text{Re }f(t + i)^2}{\cosh\pi t}dt.\endgather
$$
We need to choose $f$ so that $I\approx J\approx -K$. No matter what $f$, clearly $|K|\le J$; to achieve $J\approx -K$, the integrands in 
$J$ and $K$ must be negatives of each other, at least approximately and for most $t\in\Bbb R$ that make the integrands large. This means that $f(t+i)$ must be close to imaginary. If also $|f(t)|\approx |f(t+i)|$, then $I\approx J$. Now $f(s)=e^{\pi s/2}-e^{\pi i/4}$ satisfies both conditions and vanishes at $i/2$, but it is unbounded. Instead, with a large $\lambda\in\Bbb R$ we let
$$
f_{\lambda}(s) = \frac{e^{\pi s/2} - e^{\pi i/4}}{1+e^{\pi(s-\lambda)/2}}\ .
$$

We claim that $I(f_\lambda) \sim J(f_\lambda )\sim 2 \lambda$ and $K(f_\lambda) \sim -2\lambda$ as $\lambda \to \infty$. This will be 
verified only for $J(f_\lambda)$, the other two are treated similarly. We have
$$
J(f_\lambda) = \bigg(\int_{-\infty}^0 + \int_0^\lambda + \int_\lambda^{+\infty}\bigg)\frac{|ie^{\pi t/2} - e^{\pi i/4}|^2}{|1+ie^{\pi{(t-\lambda)}/{2}}|^2\cosh(\pi t)}dt.
$$
Since in the first integral the numerator is bounded, and in the last it is $O(\cosh \pi t)$, both integrals have bounds independent of $\lambda$. After a change of variables $\tau = {t}/{\lambda}$ in the middle integral, we obtain
$$
\align
\int_0^\lambda \frac{|ie^{{\pi t}/{2}} - e^{{\pi i}/{4}}|^2}{(1+ie^{\pi({t-\lambda})/{2}})^2\cosh(\pi t)}dt&=\lambda \int_0^1\frac{|ie^{{\pi \lambda \tau}/{2}} - e^{{\pi i}/{4}}|^2}{|1+ie^{{\pi \lambda(\tau-1)}/{2}}|^2\cosh(\pi \lambda \tau)}d\tau\\ &=2\lambda \int_0^1\frac{|i - e^{\pi({ i}/{4}-{\lambda \tau}/{2})}|^2}{|1+ie^{{\pi \lambda(\tau-1)}/{2}}|^2(1 + e^{-2\pi \lambda \tau})}d\tau.
\endalign
$$
This last expression has bounded integrand, and the dominated convergence theorem implies $J(f_\lambda) \sim 2\lambda$, as claimed. Letting $\lambda \to \infty$ in (3.5) (with $\alpha = i/2$) we obtain $0 \leq p -|q| + 2r$, and letting $p \to v_{z\overline{z}}, r \to \omega_{11}$, (3.4) follows.

To prove the sharpness of estimate (3.3), suppose that $V \subset \Bbb C$ is the unit disc and $\omega = i \partial \overline{\partial} |z|^2$. Let
$$
u(s,z) = -\frac{2\ \Im s}{\varepsilon + \Im  s}( \Re z)^2
$$
for some $\varepsilon > 0$. Clearly $u$ is bounded and continuous on $\overline{S} \times V$, $u(s,z)=0$ for all 
$(s,z) \in \Bbb R \times V$, and
$u(s,0)=0$ for all $s \in \overline{S}$. One verifies that $u$ is $\omega$-plurisubharmonic in $S \times V$ by checking that
$$
\vmatrix
 u_{s\overline{s}} & u_{s\overline{z}} \\
 u_{\overline{s}z} & 1 + u_{z\overline{z}}
\endvmatrix=0$$
and observing that $1 + u_{z\overline{z}} > 0$. This confirms that the Levi form of $|z|^2 + u$ is semi-positive everywhere. If $\varepsilon \to 0$ then $2 + v_{z\overline{z}}(0) = 2 + u_{z\overline{z}}(i,0) \to 1$ and $|v_{zz}(0)| = |u_{zz}(i,0)| \to 1$; hence the estimate (3.3) is indeed sharp.
\enddemo

\demo{Proof of Theorem 1.1}Given a $g$--invariant $v\in\cH$, suppose (1.1) has
an $\omega$--pluri\-sub\-harmonic solution $u\in C^{\pop}(\overline S\times X)$.
By Corollary 2.4, $u(s,x)=u(s,g(x))$; since $dv(x_0)=0$ is automatic for
$g$--invariant $v$, by Proposition 3.1 and by Lemma 3.2  $v$
then satisfies (3.3).
Conversely, if a $g$--invariant $v\in\cH$ does not satisfy (3.3), then (1.1) will have no $\omega$--plurisubharmonic solution $u\in C^{\pop} (\overline S\times X)$.
Such $v$ certainly exist (and form an open set among $g$--invariant potentials
in $\cH$), because the matrices $(v_{z_j\bar z_k} (x_0))=(p_{jk})$ and $(v_{z_j z_k} (x_0))=(q_{jk})$ can be arbitrarily prescribed for $g$--invariant $v\in\cH$,
as long as $(\omega_{jk}(x_0)+p_{jk})$ is positive definite, see
\cite{LV, Lemma 3.3}.
\enddemo

\enddocument
\bye